\documentclass{amsart}
\usepackage{amsmath, amssymb, latexsym, amsthm, mathrsfs}
\usepackage[all]{xy}

      \newenvironment{changemargin}[2]{\begin{list}{}{
         \setlength{\topsep}{0pt}\setlength{\leftmargin}{0pt}
         \setlength{\rightmargin}{0pt}
         \setlength{\listparindent}{\parindent}
         \setlength{\itemindent}{\parindent}
         \setlength{\parsep}{0pt plus 1pt}
         \addtolength{\leftmargin}{#1}\addtolength{\rightmargin}{#2}
         }\item }{\end{list}}

\newcommand{\Ax}{\mathsf{Ax}}

\renewcommand{\P}{\mathbb{P}}
\newcommand{\Bgood}{\B_{\mathrm{good}}}
\newcommand{\NR}{{{}^\N\R}}
\newcommand{\compactN}{\overline{\N}} 

\newcommand{\bq}{\begin{quote}}
\newcommand{\eq}{\end{quote}}

\newcommand{\CH}{the Continuum Hypothesis}

\newcommand{\inv}{^{-1}}
\newcommand{\Cantor}{{{}^\N\Z_2}}

\newcommand{\scrA}{\mathscr{A}}
\newcommand{\scrB}{\mathscr{B}}

\newcommand{\Increasing}{{{}^{\N\nearrow}\compactN}} 
\newcommand{\NcompactN}{{{}^\N\compactN}}
\newcommand{\seq}[1]{\{#1\}_{n\in\N}}

\newcommand{\cN}{\mathcal{N}}

\newcommand{\B}{\mathcal{B}}

\newcommand{\BG}{\B_\Gamma}

\newcommand{\BO}{\B_\Omega}

\newcommand{\cF}{\mathcal{F}}

\newcommand{\M}{\mathcal{M}}
\newcommand{\N}{\naturals}
\newcommand{\NN}{{{}^{\naturals}\naturals}}

\renewcommand{\O}{\mathcal{O}}
\newcommand{\Onbd}{\O_{\mathsf{nbd}}}
\newcommand{\cP}{\mathcal{P}}
\newcommand{\Q}{\rationals}
\newcommand{\R}{\reals}

\newcommand{\cU}{\mathcal{U}}
\newcommand{\Union}{\bigcup}
\newcommand{\cV}{\mathcal{V}}

\newcommand{\Z}{{\mathbb Z}}

\long\def\forget#1\forgotten{}

\renewcommand{\b}{\mathfrak{b}}

\renewcommand{\c}{\mathfrak{c}}
\renewcommand{\d}{\mathfrak{d}}

\renewcommand{\i}{\item}

\newcommand{\p}{\mathfrak{p}}

\newcommand{\w}{\omega}
\newcommand{\x}{\times}

\newcommand{\nin}{\not\in}

\newcommand{\sbst}{\subseteq}

\newcommand{\sm}{\setminus}

\renewcommand{\(}{\left(}
\renewcommand{\)}{\right)}

\renewcommand{\>}{\rangle}
\newcommand{\E}{\exists}

\newcommand{\cov}{\mathsf{cov}}

\newcommand{\cof}{\mathsf{cof}}
\newcommand{\cf}{\mathsf{cf}}

\newcommand{\impl}{\to}

\newtheorem{thm}{Theorem}

\newtheorem{fact}[thm]{Fact}
\newtheorem{prob}[thm]{Problem}
\newtheorem{lem}[thm]{Lemma}
\newtheorem{cor}[thm]{Corollary}

\theoremstyle{definition}

\theoremstyle{remark}
\newtheorem{rem}[thm]{Remark}
\newcommand{\be}{\begin{enumerate}}
\newcommand{\ee}{\end{enumerate}}
\newcommand{\bi}{\begin{itemize}}
\newcommand{\ei}{\end{itemize}}

\forget
\forgotten

\newcommand{\sone}{\mathsf{S}_1}
\newcommand{\sfin}{\mathsf{S}_{fin}}
\newcommand{\ufin}{\mathsf{U}_{fin}}

    \newcommand{\gfin}{\mathsf{G}_{fin}}
\newcommand{\naturals}{{\mathbb N}}
\newcommand{\reals}{{\mathbb R}}
\newcommand{\rationals}{{\mathbb Q}}

\author{Boaz Tsaban}
\address{Department of Applied Mathematics and Computer Science,
Weizmann Institute of Science,
Rehovot 76100,
Israel}
\email{boaz.tsaban@weizmann.ac.il}
\urladdr{http://www.cs.biu.ac.il/\~{}tsaban}

\title[$o$-bounded and combinatorial groups]{%
$o$-bounded groups and other topological groups with strong combinatorial properties}

\begin{document}

\begin{abstract}
We construct several topological groups with very strong combinatorial properties.
In particular, we give simple examples of subgroups of $\R$ (thus strictly $o$-bounded)
which have the Menger and Hurewicz properties but are not $\sigma$-compact,
and show that the product of two
$o$-bounded subgroups of $\NR$ may fail to be $o$-bounded,
even when they satisfy the stronger property $\sone(\BO,\BO)$.
This solves a problem of Tka\v{c}enko and Hernandez,
and extends independent solutions of Krawczyk and Michalewski
and of Banakh, Nickolas, and Sanchis.
We also construct separable metrizable groups $G$ of size continuum
such that every countable Borel
$\w$-cover of $G$ contains a $\gamma$-cover of $G$.
\end{abstract}


\maketitle

\vspace{-0.6cm}

\section{Introduction}

In \cite{coc1, coc2}, a unified framework for topological diagonalizations
is established, which turns out closely related to several notions which
appear in a recently flourishing study of topological groups in terms of
their covering properties (see, e.g., \cite{TkaIntro, Hernandez, HRT, KMexample, KMlinear}
and references therein).
A comprehensive study of these interrelations is currently being carried by
Babinkostova, Ko\v{c}inac, and Scheepers \cite{coc11}.
The purpose of this paper is to adopt several recent construction
techniques from the general theory of topological diagonalizations
to the theory of topological groups.

\subsection{Topological diagonalizations}
We briefly describe the general framework.
Let $X$ be a topological space.
An open cover $\cU$ of $X$ is
an \emph{$\omega$-cover} of $X$ if $X$ is not in $\cU$ and for
       each finite subset $F$ of $X$, there is
       a set $U\in\cU$ such that $F\subseteq U$.
$\cU$ is a \emph{$\gamma$-cover} of $X$ if it is infinite and for each $x$ in
       $X$, $x\in U$ for all but finitely many $U\in\cU$.
Let $\O$, $\Omega$, and $\Gamma$ denote the collections of all countable open
covers, $\omega$-covers, and $\gamma$-covers of $X$, respectively.
Let $\scrA$ and $\scrB$ be collections of covers of a space $X$.
Following are selection hypotheses which $X$ might satisfy or not
satisfy.
\begin{itemize}
\item[$\sone(\scrA,\scrB)$:]{
For each sequence $\seq{\cU_n}$ of members of $\scrA$,
there exist members $U_n\in\cU_n$, $n\in\N$, such that $\seq{U_n}\in\scrB$.
}
\item[$\sfin(\scrA,\scrB)$:]{
For each sequence $\seq{\cU_n}$
of members of $\scrA$, there exist finite (possibly empty)
subsets $\cF_n\sbst\cU_n$, $n\in\N$, such that $\Union_{n\in\N}\cF_n\in\scrB$.
}
\item[$\ufin(\scrA,\scrB)$:]{
For each sequence $\seq{\cU_n}$ of members of $\scrA$
which do not contain a finite subcover,
there exist finite (possibly empty) subsets $\cF_n\sbst\cU_n$, $n\in\N$,
such that $\seq{\cup\cF_n}\in\scrB$.
}
\end{itemize}

\begin{figure}[!ht]
{
\begin{changemargin}{-2cm}{-1cm}
\begin{center}
$\xymatrix@R=10pt{
&
&
& \ufin(\O,\Gamma)\ar[r]
& \ufin(\O,\Omega)\ar[rr]
& & \sfin(\O,\O)
\\
&
&
& \sfin(\Gamma,\Omega)\ar[ur]
\\
& \sone(\Gamma,\Gamma)\ar[r]\ar[uurr]
& \sone(\Gamma,\Omega)\ar[rr]\ar[ur]
& & \sone(\Gamma,\O)\ar[uurr]
\\
&
&
& \sfin(\Omega,\Omega)\ar'[u][uu]
\\
& \sone(\Omega,\Gamma)\ar[r]\ar[uu]
& \sone(\Omega,\Omega)\ar[uu]\ar[rr]\ar[ur]
& & \sone(\O,\O)\ar[uu]
}$
\caption{The Scheepers Diagram}\label{SchDiagram}
\end{center}
\end{changemargin}
}
\end{figure}
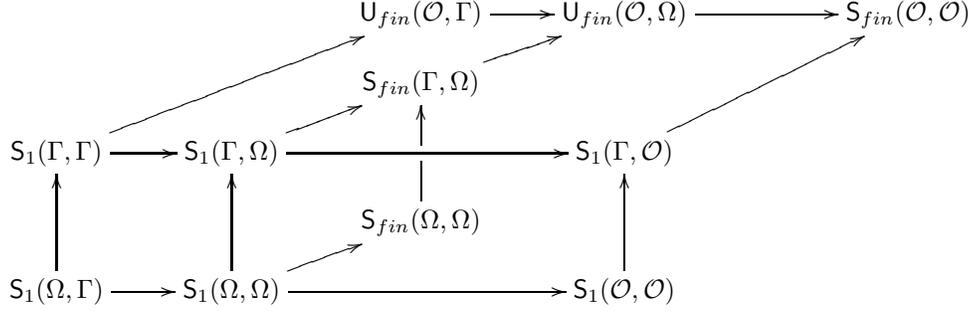

Some of the properties defined in this manner
were studied earlier by Hurewicz ($\ufin(\O,\Gamma)$), Menger ($\sfin(\O, \O)$),
Rothberger ($\sone(\O, \O)$, traditionally known as the $C''$ property),
Gerlits and Nagy ($\sone(\Omega,\Gamma)$, traditionally known as the $\gamma$-property),
and others.
Many equivalences hold among these properties, and the surviving ones
appear in Figure \ref{SchDiagram} (where an arrow denotes implication),
to which no arrow can be added except perhaps from
$\ufin(\O,\Gamma)$ or $\ufin(\O,\Omega)$ to $\sfin(\Gamma,\Omega)$
\cite{coc2}.

Each selection principle has a naturally associated game,
but we will restrict attention to the game
$\gfin(\scrA,\scrB)$, which is played as follows:
In the $n$th inning, ONE chooses an element
$\cU_n$ of $\scrA$ and then TWO responds by choosing
a finite subset $\cF_n$ of $\cU_n$.
They play an inning per natural number.
A play $(\cU_0, \cF_0, \cU_1, \cF_1\dots)$ is won by TWO if
$\Union_n\cF_n\in\scrB$; otherwise ONE wins.
We will write $\mbox{ONE}\uparrow\gfin(\scrA,\scrB)$ (respectively, $\mbox{TWO}\uparrow\gfin(\scrA,\scrB)$) for
``ONE (respectively, TWO) has a winning strategy in $\gfin(\scrA,\scrB)$''.

\subsection{$o$-bounded groups}
Okunev introduced the following notion as an approximation of $\sigma$-compact
groups: A topological group $G$ is \emph{$o$-bounded} if for each sequence
$\seq{U_n}$ of neighborhoods of the unit element of $G$, there exists
a sequence $\seq{F_n}$ of finite subsets of $G$ such that
$G=\Union_n F_n\cdot U_n$.
It is possible to state this definition in the language of selection principles.
Let $\Onbd$ denote the covers of $G$ of the form $\{g\cdot U : g\in G\}$
where $U$ is a neighborhood of the unit element of $G$.
Then $G$ is $o$-bounded if, and only if, $G$ satisfies $\sfin(\Onbd,\O)$.

According to Tka\v{c}enko, a topological group $G$ is \emph{strictly $o$-bounded}
if TWO has a winning strategy in the game $\gfin(\Onbd,\O)$.
Clearly, every subgroup of a $\sigma$-compact group is strictly
$o$-bounded, but the converse does not hold \cite{Hernandez}.

\subsubsection*{Notational convention}
For sets $X,Y$, ${{}^XY}$ denotes the collection of all functions from $X$ to $Y$.
If $Y$ is a topological space, then the topology on ${{}^XY}$ is the Tychonoff product
topology.

\section{Two almost $\sigma$-compact subgroups of $\R$}
The \emph{Baire space} $\NN$ 
is (quasi)ordered by eventual dominance: $f\le^* g$ if
$f(n)\le g(n)$ for all but finitely many $n$.
A subset of $\NN$ is \emph{dominating} if it is cofinal in $\NN$ with respect to $\le^*$.
If a subset of $\NN$ is unbounded with respect to $\le^*$ then we simply say that it is
\emph{unbounded}.
Let $\b$ (respectively, $\d$) denote the minimal cardinality of an unbounded
(respectively, dominating) subset of $\NN$.

We use the following setting from \cite{ideals}.
Let $\compactN=\N\cup\{\infty\}$ be the one point compactification of $\N$. (A subset
$A\sbst\compactN$ is open if: $A\sbst\N$, or $\infty\in A$ and $A$
is cofinite.) Let $\Increasing \sbst \NcompactN$ consist of the
\emph{nondecreasing} functions $f\in\NcompactN$.
$\Increasing$ is homeomorphic to the Cantor set of reals.
For each increasing finite sequence
$s$ of natural numbers, let $q_{s} \in \Increasing$ be defined by
$q_{s}(k) = s(k)$ if $k <|s|$, and $q_{s}(k) = \infty$ otherwise.
Let $Q$ be the collection of all these elements $q_s$.

\begin{thm}\label{ZFCgp}
There exists a non $\sigma$-compact subgroup $G_H$ of $\R$ of
cardinality $\b$ such that all finite powers of $G_H$
satisfy $\ufin(\O,\Gamma)$ (in particular, they satisfy $\sfin(\Omega,\Omega)$).
\end{thm}
\begin{proof}
Let $B=\{f_\alpha : \alpha < \b\}\sbst\NN$ be a $\le^*$-unbounded
set of strictly increasing elements of $\NN$ which forms a \emph{$\b$-scale}
(that is, for each $\alpha<\beta$, $f_\alpha \leq^* f_\beta$),
and set $H = B\cup Q$.
In \cite{ideals} it is proved that
all finite powers of $H$ satisfy $\ufin(\O,\Gamma)$.

Think of $H$ as a set of real numbers.
For each $n$, the set
$$G_n = \{\alpha_1g_1+\dots+\alpha_ng_n : \alpha_1,\dots,\alpha_n\in\Z,\ g_1,\dots,g_n\in H\}$$
is a union of countably many continuous images of $H^n$, thus for each $k$, $(G_n)^k$
is a union of countably many continuous images of $H^{nk}$.
As the property $\ufin(\O,\Gamma)$ is preserved under taking continuous images
and countable unions \cite{coc2, AddQuad}, we have that each set $(G_n)^k$ satisfies
$\ufin(\O,\Gamma)$.

Take $G_H=\<H\>$.
Then $G_H=\Union_{n\in\N}G_n$ and for each $n$, $G_n\sbst G_{n+1}$.
Thus, $(G_H)^k = \Union_{n\in\N} (G_n)^k$ for each $k$;
therefore $(G_H)^k$ satisfies $\ufin(\O,\Gamma)$ for each $k$.
Now, satisfying $\ufin(\O,\O)$ in all finite powers
implies $\sfin(\Omega,\Omega)$ \cite{coc2}.

It was observed by Pol and Zdomskyy, that one can make sure that
$G_H$ is not $\sigma$-compact by embedding $\Increasing$ in a Cantor set of reals, $C$,
that is \emph{linearly independent} over $\Q$. This way, $H=\<H\>\cap C$ is a closed
subset of $\<H\>=G_H$. Since $H$ is not $\sigma$-compact \cite{ideals}, we get that $G_H$
cannot be $\sigma$-compact.
\end{proof}

\begin{prob}
Does $G_H$ satisfy $\sone(\Gamma,\Gamma)$?
\end{prob}

\begin{cor}
$\mbox{TWO}\uparrow\gfin(\O,\O)$ is strictly stronger than
$\mbox{TWO}\uparrow\gfin(\Onbd,\O)$ (strict $o$-boundedness).
\end{cor}
\begin{proof}
By a well known theorem of Telg\'arsky, $\mbox{TWO}\uparrow\gfin(\O,\O)$
if, and only if, the space $G$ is $\sigma$-compact.
Contrast this with Theorem \ref{ZFCgp}.
\end{proof}

Theorem \ref{ZFCgp} has a group-theoretic consequence.

\begin{cor}
Assume that for each $m$, $\seq{g_{m,n}+U_m}$ is a cover of $G_H$.
Then there exists a sequence
$\seq{m_n}$ such that
$$G_H=\Union_k\bigcap_{n>k}(\{g_{n,1},\dots,g_{n,m_n}\} + U_n).$$
In fact, this property is satisfied by all finite powers of $G_H$.
\end{cor}

Let $D =\{g_\alpha : \alpha < \d\}$ be a dominating subset
of $\NN$ where each $g_\alpha$ is increasing, and
for each $f\in\NN$ there exists $\alpha_0<\d$ such that for any
finite set $F\sbst\d\setminus\alpha_0$,
$f(n) < \min\{g_\beta(n) : \beta\in F\}$
for infinitely many $n$.
Such a set was constructed in \cite{ideals}.

A subset $F$ of $\NN$ is \emph{finitely-dominating}
if for each $g\in\NN$ there exist $k$ and $f_1,\dots,f_k\in\NN$
such that $g(n)\le^* \max\{f_1(n),\dots,f_k(n)\}$.
For conciseness, we use the following shortened notation.
\bi
\i[$\Ax$:]
Either a union of less than $\d$ many not dominating sets is not
dominating (in other words, $\b=\d$), or else a union of less
than $\d$ many not finitely-dominating sets
of increasing functions is not finitely-dominating.
\ei
In \cite{ideals} it is shown that
$\Ax$ implies that $M=D\cup Q$ satisfies $\sfin(\Omega,\Omega)$.
(Observe that $\sfin(\Omega,\Omega)$ is preserved under taking finite powers \cite{coc2}.)
Assuming $\Ax$, one shows as in Theorem \ref{ZFCgp} that all finite powers of $G_M=\<M\>$
satisfy $\ufin(\O,\O)$ and gets the following.
\begin{thm}\label{mengrp}
Assume that $\Ax$ holds. Then:
\be
\i There exists a non $\sigma$-compact subgroup $G_M$ of $\R$
of cardinality $\d$ such that $G_M$ satisfies $\sfin(\Omega,\Omega)$.
\i Assume that for each $m$, $\seq{g_{m,n}+U_m}$ is a cover of $G_M$.
Then there exists a sequence
$\seq{m_n}$ such that for each finite subset $F$ of $G_M$, there exists
$n$ such that $F\sbst \{g_{n,1},\dots,g_{n,m_n}\} + U_n$.
Moreover, this holds for all finite powers of $G_M$.
\ee
\end{thm}

It follows from the next section that the hypothesis $\Ax$ is not necessary to prove
Theorem \ref{mengrp} (namely, it also follows from the incomparable assumption $\cov(\M)=\c$).

\begin{prob}
Is Theorem \ref{mengrp} provable in ZFC?
\end{prob}

\section{Products of $o$-bounded groups}

In Problem 3.2 of \cite{TkaIntro} and Problem 5.2 of \cite{Hernandez}
it is asked whether the (Tychonoff) product of two $o$-bounded groups is $o$-bounded.
We give a negative answer.
A negative answer was independently given by Krawczyk and Michalewski \cite{KMlinear},
but our result is stronger in the following sense:
Let $G$ be a topological group such that all finite powers of $G$ are Lindel\"of.
Then each open $\w$-cover of $G$ contains a \emph{countable} $\w$-cover of $G$.
Let $\BO$ denote the collection of all \emph{countable Borel} $\w$-covers of $G$.
In this case,
$$\sone(\BO,\BO) \impl \sone(\Omega,\Omega)\impl \sfin(\O,\O)$$
where no implication can be reversed \cite{coc2, CBC}, and the last property
(the Menger property) implies $\sfin(\Onbd,\O)$ ($o$-boundedness).
In \cite{KMlinear} it is proved that, assuming $\cov(\M)=\c$ (this is a small portion of
\CH), there exist groups $G_1$ and $G_2$ satisfying the Menger property $\sfin(\O,\O)$,
such that $G_1\x G_2$ is not $o$-bounded.
We use the same hypothesis to show that there exist such
groups satisfying $\sone(\BO,\BO)$.

\begin{fact}
Assume that $G$ is a $\le^*$-dominating subgroup of $\NR$.
Then $G$ is not $o$-bounded.
\end{fact}
\begin{proof}
As $\N$ can be partitioned into infinitely many infinite sets,
the following holds.
\begin{lem}
For each $o$-bounded group $G$ and sequence
$\seq{U_n}$ of neighborhoods of the identity of $G$,
there exists a sequence $\seq{F_n}$ of finite subsets of $G$
with $G=\bigcap_m\Union_{n>m}F_n\cdot G$.
\end{lem}
Consider the open sets $U_n = \{f\in\NR : |f(n)|<1\}\sbst\NR$.
For each sequence $\seq{F_n}$ of finite subsets of
$\NR$, $\Union_nF_n+U_n$ is $\le^*$-bounded in $\NR$.
Let $h\in\NR$ be a witness for that.
As $G$ is dominating, there exists $g\in G$ such that $h\le^* g$.
Then $g\nin\bigcap_m\Union_{n>m}F_n\cdot G$.
\end{proof}

Let $\c=|\R|$. The assertion $\cov(\M)=\c$ means that $\R$ (or any complete, separable,
metric space) is not the union of less than continuum many meager (=first category) sets.
We say that $L$ is a \emph{$\kappa$-Luzin group} if $|L|\ge\kappa$, and
for each meager set $M$ in $L$, $|L\cap M|<\kappa$.
\begin{thm}\label{prods}
Assume that $\cov(\M)=\c$.
Then there exist $\c$-Luzin subgroups $L_1,L_2$ of $\NR$ satisfying $\sone(\BO,\BO)$,
such that the group $L_1\x L_2$ is not $o$-bounded.
(These Luzin groups are, in fact, linear vector spaces over $\Q$.)
\end{thm}
\begin{proof}
We extend the technique of \cite{CBC, AddQuad}.
We stress that there exists a much easier proof if we only require
that $L_1$ and $L_2$ satisfy $\sone(\Omega,\Omega)$; however we
do not supply this easier proof to avoid repetitions.

A cover $\cU$ of $X$ is \emph{$\w$-fat} if for each finite $F\sbst X$
and each finite family $\cF$ of
nonempty open sets, there exists $U\in\cU$
such that $F\sbst U$ and for each $O\in\cF$, $U\cap O$ is not meager.
Let $\M$ denote the meager subsets of $\NR$.
\begin{lem}[\cite{huremen2}]\label{fatlemma}
Assume that $\cU$ is an $\w$-fat cover of a set $X\sbst\NR$.
Then:
\be
\i $\cup\cU$ is comeager,
\i For each finite $F\sbst G$ and finite family $\cF$ of
nonempty open sets,
$$\cU_{F,\cF} := \{U\in\cU : F\sbst U\mbox{ and for each $O\in\cF$, }U\cap O\nin\M\}$$
is an $\w$-fat cover of $G$.
Consequently, $\cup\cU_{F,\cF}$ is comeager.
\ee
\end{lem}

For shortness, we will say that a cover $\cU$ of $X$ is \emph{good} if:
For each finite $A\sbst\Q\sm\{0\}$ and finite $B\sbst G$,
the family
$$\cU^{A,B} := \left\{\bigcap_{q\in A, g\in B} q(U-g) : U\in\cU\right\}$$
is an $\w$-fat cover of $X$.
These covers allow the inductive construction hinted in the following lemma.

\begin{lem}\label{addelement}
Assume that $\cU$ is a good cover of a group $G\sbst\NR$.
Then for each element $x$ in the intersection of all sets of the form
$\cup(\cU^{A,B})_{F,\cF}$
where the members of $\cF$ are \emph{basic} open sets,
$\cU$ is a good cover of the group $G+\Q x$.
\end{lem}
\begin{proof}
Fix finite sets $A\sbst\Q\sm\{0\}$ and $B\sbst G$.
We may assume that $1\in A$.
We must show that $\cU^{A,B}$ is an $\w$-fat cover of $G+\Q x$.
Let $F$ be a finite subset of $G+\Q x$, and $\cF$ be a finite family of
nonempty open sets. By moving to subsets we may assume that all members of
$\cF$ are basic open sets.

Choose finite sets $\tilde A\sbst\Q\sm\{0\}$, $\tilde B\sbst G$, and $\tilde F\sbst G$,
such that $F\sbst (\tilde B + \tilde A x)\cup \tilde F$, $1\in \tilde A$, and $0\in\tilde B$.
As $x\in\cup(\cU^{\tilde A\inv A, B+A\inv\tilde B})_{\tilde F,\cF}$,
there exists $U\in\cU$ such that
$x\in \tilde V:=\bigcap_{q\in \tilde A\inv A, g\in B+A\inv\tilde B}q(U-g)$,
$\tilde F\sbst\tilde V$, and for each $O\in\cF$, $\tilde V\cap O$ is not meager.
Take $V=\bigcap_{q\in A, g\in B}q(U-g)$. Then $\tilde V\sbst V$, thus
$\tilde F\sbst V$ and for each $O\in\cF$, $V\cap O$ is not meager.
Now, $x\in\tilde V$, thus for each $\tilde a\in\tilde A$ and $\tilde b\in \tilde B$,
$x\in\bigcap_{q\in \tilde a\inv A, g\in B+A\inv\tilde b}q(U-g)$, thus
$x\in\bigcap_{q\in A, g\in B}\tilde a\inv q(U-(g+q\inv\tilde b))$, therefore
$x\in\bigcap_{q\in A, g\in B}\tilde a\inv (q(U-g)-\tilde b)$, thus
$\tilde a x+\tilde b\in V$. This shows that
$F\sbst (\tilde A x + \tilde B)\cup\tilde F\sbst V\in\cU^{A,B}$,
and we are done.
\end{proof}

Since we are going to construct Luzin groups, the following lemma
tells that we need not consider covers which are not good.
\begin{lem}\label{denselusin}
Assume that $L$ is a subgroup of $\NR$ such that $\Q\cdot L\sbst L$ and
for each nonempty basic open set $O$, $L\cap O$ is not meager.
Then every countable Borel $\w$-cover $\cU$ of $L$ is a good cover of $L$.
\end{lem}
\begin{proof}
Assume that $\cU$ is a countable collection of Borel sets which is
not a good cover of $L$.
Then there exist finite sets $A\sbst\Q\sm\{0\}$, $B\sbst L$,
$F\sbst L$, and $\cF$ of nonempty open sets
such that for each $V\in\cU^{A,B}$ containing $F$, $V\cap O$ is meager for some $O\in\cF$.
For each $O\in\cF$ let
$$M_O = \cup\{V\in\cU^{A,B} : F\sbst V\mbox{ and }V\cap O\in\M\}.$$
Then $M_O\cap O$ is meager,
thus there exists $x_O\in (L\cap O)\sm M_O$.
Then $F\cup\{x_O : O\in\cF\}$ is not covered by any $V\in\cU^{A,B}$.
We will show that this cannot be the case.

Put $F' = A\inv F + B$.
As $\Q\cdot L\sbst L$, $F'$ is a finite subset of $L$.
Thus, there exists $U\in\cU$ such that $F'\sbst U$.
Consequently, for each $q\in A$ and $g\in B$,
$x\in q(U-g)$ for each $x\in F$, that is,
$F\sbst \bigcap_{q\in A, g\in B} q(U-g)\in\cU^{A,B}$.
\end{proof}

We need one more lemma. Denote the collection of countable Borel
good covers of a set $X$ by $\Bgood$.
\begin{lem}\label{s1goodgood}
If $|X|<\cov(\M)$, then $X$ satisfies $\sone(\Bgood,\Bgood)$.
\end{lem}
\begin{proof}
Assume that $|X|<\cov(\M)$, and let $\seq{\cU_n}$ be a sequence of
countable Borel good covers of $X$.
Enumerate each cover $\cU_n$ by $\{U^n_k\}_{k\in\N}$.
Let $\seq{Y_n}$ be a partition of $\N$ into infinitely many
infinite sets. For each $m$, let $y_m\in\NN$ be an increasing
enumeration of $Y_m$.
Let $\seq{\cF_n}$ be an enumeration of all finite families of
nonempty basic open sets.

For finite sets $F,B\sbst X$ and $A\sbst\Q\sm\{0\}$,
and each $m$ define a function
$\Psi^{A,B}_{F,m}\in\NN$ by
$$\Psi^{A,B}_{F,m}(n)=\min\{k : F\sbst V:= \bigcap_{q\in A,g\in B}q(U^{y_m(n)}_k-g)\mbox{ and }(\forall O\in\cF_m)\
V\cap O\nin\M\}.$$
Since there are less than $\cov(\M)$ many functions $\Psi^{A,B}_{F,m}$,
there exists by \cite{covM} a function $f\in\NN$ such that
for each $m$, $F$, $A$, and $B$, $\Psi^{A,B}_{F,m}(n)=f(n)$ for infinitely many $n$.
Consequently, $\cV=\{U^{y_m(n)}_{f(n)} : m,n\in\N\}$ is a good cover of $X$.
\end{proof}

We are now ready to carry out the construction.
Let $\NR=\{y_{\alpha}:\alpha<\c\}$,
$\{M_{\alpha}:\alpha<\c\}$ be all $F_\sigma$ meager subsets of $\NR$, and
$\{\seq{\cU^\alpha_n} : \alpha<\c\}$ be all sequences of countable families of
Borel sets.
Let $\{O_k : k\in\N\}$ and $\{\cF_m: m\in\N\}$ be all nonempty basic open sets
and all finite families of nonempty basic open sets, respectively, in $\NR$.

We construct $L_1$ and $L_2$ by induction on $\alpha<\c$ as follows.
At stage $\alpha \ge 0$ set
$L^i_\alpha  = \Union_{\beta<\alpha}L^i_\beta$
and consider the sequence $\seq{\cU^\alpha_n}$.
Say that $\alpha$ is $i$-good if for each $n$
$\cU^\alpha_n$ is a good cover of $L^i_\alpha$.
In this case, by Lemma \ref{s1goodgood} there exist elements
$U^{\alpha,i}_n\in\cU^\alpha_n$ such that $\cU^{\alpha,i} = \seq{U^{\alpha,i}_n}$ is
a good cover of $L^i_\alpha$.
We make the inductive hypothesis that
for each $i$-good $\beta<\alpha$,
$\cU^{\beta,i}$ is a good cover of $L^i_\alpha$.
For finite sets $F,B\sbst L^i_\alpha$ and $A\sbst\Q\sm\{0\}$,
each $i$-good $\beta\le\alpha$, and each $m$ define
$$G^{\beta,i}_{A,B,F,m}=\cup((\cU^{\beta,i})^{A,B})_{F,\cF_m}.$$
As $\cU^{\beta,i}$ is a good cover of $X^i_\alpha$, $(\cU^{\beta,i})^{A,B}$
is $\w$-fat cover of $L^i_\alpha$, and
by Lemma \ref{fatlemma}, $G^{\beta,i}_{A,B,F,m}$ is comeager in $\NR$.
Set
$$Y_\alpha=\Union_{\beta<\alpha}M_\beta\cup
\Union\left\{\NR\sm G^{\beta,i}_{A,B,F,m} :
\txt{
$i<2,\mbox{ $i$-good }\beta\le\alpha, m\in\N,$\\
finite sets $F,B\sbst L^i_\alpha, A\sbst\Q\sm\{0\}$
}
\right\},$$
and $Y_\alpha^* = \{x\in\NR : (\E y\in
Y_\alpha)\ x=^* y\}$ (where $x =^* y$ means that $x(n)=y(n)$ for
all but finitely many $n$.) Then $Y_\alpha^*$ is a union of less
than $\cov(\M)$ many meager sets.
\begin{lem}[\cite{huremen2}]\label{x+y=z}
If $X$ is a union of less than $\cov(\M)$ many
meager sets in $\NR$,
then for each $x\in\NR$ there exist $y,z\in\NR\sm X$ such that $y+z=x$.
\end{lem}
Use Lemma \ref{x+y=z} to pick
$x^0_\alpha,x^1_\alpha\in\NR\sm Y_\alpha^*$ such that
$x^0_\alpha+x^1_\alpha=y_\alpha$. Let $k = \alpha \bmod \omega$,
and change a finite initial segment of $x^0_\alpha$ and
$x^1_\alpha$ so that they both become members of $O_k$. Then
$x^0_\alpha,x^1_\alpha\in O_k\sm Y_\alpha$, and
$x^0_\alpha+x^1_\alpha =^* y_\alpha$.
Finally, define $X^i_{\alpha+1} = L^i_\alpha + \Q\cdot x^i_\alpha$.
By Lemma \ref{addelement}, the inductive hypothesis is preserved.
This completes the construction.

Take $L_i = \Union_{\alpha<\c}L^i_\alpha$, $i=1,2$.
By Lemma \ref{s1goodgood},
each $L_i$ satisfies $\sone(\Bgood,\Bgood)$,
and by the construction,
its intersection with each nonempty basic open set has size $\c$.
By Lemma \ref{denselusin}, $\Bgood = \BO$ for $L_i$.
Finally, $L_1+L_2$ (a homomorphic image of $L_1\x L_2$)
is a dominating subset of $\NR$, thus $L_1\x L_2$ is not $o$-bounded.
\end{proof}

\begin{rem}
None of the $o$-bounded groups $L_1$ and $L_2$ in Theorem \ref{prods} is
strictly $o$-bounded. By an unpublished result of Michalewski,
every metrizable strictly $o$-bounded group is a subgroup of a $\sigma$-compact
group, and in \cite[Theorem 5.3]{Hernandez} it is proved that a product of
such a group with an $o$-bounded group is $o$-bounded.
\end{rem}

The following consequence of Theorem \ref{prods} seems nontrivial.
Say that a subset $S$ of a topological group $G$ is
$\aleph_0$-bounding if there exists a countable
set $F\sbst G$ such that $F\cdot S=G$.
For example, $G$ is \emph{$\aleph_0$-bounded}
if each nonempty open set in $G$
is $\aleph_0$-bounding.
The first property in the following corollary may be called
\emph{Borel $o$-boundedness}.
This property is more interesting when the group in question
is $\aleph_0$-bounded, in which case it is stronger than
$o$-boundedness.

\begin{cor}
Assume that $\cov(\M)=\c$.
Then there exists an $\aleph_0$-bounded topological group $G$
of size continuum such that:
\be
\i For each sequence $\seq{B_n}$ of
$\aleph_0$-bounding Borel sets in $G$,
there exits a sequence $\seq{F_n}$ of finite
subsets of $G$ such that $G=\Union_n F_n\cdot B_n$.
\i Moreover, the sequence in (1) will have the
property that for each finite $F\sbst G$ there exists
$n$ such that $F\sbst F_n\cdot U_n$.
\ee
\end{cor}

\begin{rem}
Banakh, Nickolas, and Sanchis have also, independently, proved the consistency of
$o$-bounded groups not being closed under taking finite prodocts
(however, they do not consider stronger combinatorial properties as done in
\cite{KMlinear} and here).
Their construction uses ultrafilters which are not nearly coherent --
see \cite{BNS}.
\end{rem}

\section{Groups satisfying $\sone(\BO,\BG)$ or $\sone(\BG,\BG)$}

One may wonder whether Theorem \ref{prods} can be strengthened further
so that $L_1$ and $L_2$ will satisfy a stronger property.
By inspection of the Scheepers Diagram \ref{SchDiagram}, the only
candidate for a stronger property (among the ones considered in
this paper) is $\sone(\BO,\BG)$. This is far from possible:
A result of \cite{coc11} (see \cite[Theorem 32]{LecceSurvey})
implies that whenever $G_1$ is a topological group satisfying $\sone(\Omega,\Gamma)$,
and $G_2$ is $o$-bounded, the group $G_1\x G_2$ is $o$-bounded.
However strong, though, the notion of a topological group satisfying
$\sone(\BO,\BG)$ is not trivial.
\begin{thm}
For each cardinal $\kappa$ with $\cf(\kappa)>\aleph_0$,
it is consistent that $\c=\kappa$ and there exists a topological
subgroup of $\R$ of size $\c$, which satisfies $\sone(\BO,\BG)$.
\end{thm}
\begin{proof}
This is really a theorem of Miller:
Let $M$ be a countable standard model of ZFC
satisfying $\c=\kappa$.
In \cite{MilNonGamma} it is proved that
there exists a ccc poset $\P$ in $M$ of size continuum
(so that forcing with $\P$ does not change the size of the continuum)
such that the old reals $M\cap\R$ satisfy $\sone(\BO,\BG)$ in $V^\P$.

But observe that, as the operations of addition and substraction in $\R$
are absolute, $G=M\cap\R$ is a group in $V^\P$.
\end{proof}


\begin{prob}\label{gammagp}
Does \CH{} imply the existence of a separable metrizable
group of size $\c$ which satisfies $\sone(\BO,\BG)$?
\end{prob}

We have some approximate results.
For a sequence $\seq{X_n}$ of subsets of $X$, define
$\liminf X_n = \Union_m\bigcap_{n\ge m} X_n$.
For a family $\cF$ of subsets of $X$, $L(\cF)$ denotes
its closure under the operation $\liminf$.
According to \cite{GN},
$X$ is a \emph{$\delta$-set} if
for each $\w$-cover $\cU$ of $X$, $X\in L(\cU)$.
It is easy to see that the $\delta$-property is preserved under taking countable
\emph{increasing} unions.
Clearly $\sone(\Omega,\Gamma)$ implies the $\delta$-property.
The converse is still an open problem \cite{futurespm}.
If a $\delta$-set is a group, we will call it a $\delta$-group.

Let $\Z_2$ denote the usual group $\{0,1\}$ with modulo $2$ addition.
\begin{thm}
Assume that $\p=\c$. Then there exists a subgroup $G$ of $\Cantor$
such that for each $k$ $G^k$
is a countable increasing union of sets satisfying $\sone(\Omega,\Gamma)$.
In particular, all finite powers of $G$ are $\delta$-groups.
\end{thm}
\begin{proof}
As $\p=\c$, there exists a subset $X$ of $\Cantor$ of size continuum which
satisfies $\sone(\Omega,\Gamma)$ \cite{GM}.
We may assume that $0\in X$. Consequently,
$$G := \<X\> = \Union_{n\in\N}\{x_1+\dots+x_n : x_1,\dots,x_n\in X\}$$
is a countable increasing union of continuous images of powers of
$X$. But the property $\sone(\Omega,\Gamma)$ is closed under taking
finite powers and continuous images \cite{coc2}.
Observe that each finite power of $G$ is the countable increasing union
of the same power of the original sets.
\end{proof}

Assuming \CH{}, there exists a set of reals $X$ of size $\c$
satisfying $\sone(\BO,\BG)$ (e.g., \cite{MilNonGamma}).
It is an open problem whether $\sone(\BO,\BG)$ is provably preserved under
taking finite powers \cite{futurespm}. If it is, then we have a positive
answer to Problem \ref{gammagp}.

\begin{thm}\label{borgamma}
Assume that $X\sbst\Cantor$ satisfies $\sone(\BO,\BG)$ in all finite powers.
Then $G=\<X\>$ satisfies $\sone(\BO,\BG)$ in all finite powers.
\end{thm}
\begin{proof}
By the above arguments, it suffices to prove the following.
\begin{lem}
$\sone(\BO,\BG)$ is preserved under taking countable increasing unions.
\end{lem}
To prove the lemma, assume that $X=\Union_n X_n$ is an increasing union,
and observe that $\sone(\BO,\BG)$ implies
$\sone(\BG,\BG)$, which in turn implies that all Borel images of
each $X_k$ in $\NN$ are bounded \cite{CBC}.

Assume that $\cU_n=\{U^n_m : m\in\N\}$, $n\in\N$, are countable Borel $\w$-covers of $X$.
By $\sone(\BO,\BG)$, we may assume that each $\cU_n$ is a $\gamma$-cover of $X_n$.
For each $k$, define a function $\Psi_k$ from $X_k$ to $\NN$ so that for each $x$ and $n$:
$$\Psi_k(x)(n) = \min\{i : \(\forall m\geq i\)\ x\in U^n_m\}.$$
$\Psi_k$ is a Borel map, thus $\Psi_k[X]$ is bounded.
Consequently, $\Union_k\Psi_k[X]$ is bounded, say by the sequence $m_n$.
Then $\seq{U^n_{m_n}}$ is a $\gamma$-cover of $X$, as required.
\end{proof}

We conclude the paper with a group satisfying $\sone(\BG,\BG)$ in all finite powers.
Let $\cN$ denote the collection of null (Lebesgue measure zero) sets of reals.
$\cov(\cN)$ is the minimal size of a cover of $\R$ by null sets, and
$\cof(\cN)$ is the minimal size of a cofinal family in $\cN$ with respect to inclusion.
Let $\kappa$ be an uncountable cardinal.
$S\sbst\R$ is a $\kappa$-Sierpi\'nski set if $|S|\ge\kappa$ and
for each null set $N$, $|S\cap N|<\kappa$. $\b$-Sierpi\'nski sets satisfy
$\sone(\BG,\BG)$, but powers of $\kappa$-Sierpi\'nski sets are never
$\kappa$-Sierpi\'nski sets.
\begin{thm}
Assume that $\cov(\cN)=\cof(\cN)=\b$.
Then there exists a (non $\sigma$-compact)
group $G\sbst\R$ of size $\b$, such that
all finite powers of $G$ satisfy $\sone(\BG,\BG)$.
Consequently, all finite powers of $G$ also satisfy
$\sfin(\BO,\BO)$.
\end{thm}
\begin{proof}
The proof is similar to that of Theorem 7.2 of \cite{BRR2}.
\begin{lem}\label{sier}
Assume that $\cov(\cN)=\cof(\cN)$.
Then there exists a $\cov(\cN)$-Sierpi\'nski set $S$ such that
for each $k$ and each null set $N$ in $\R^k$,
$S^k \cap N$ is contained in a union of less than $\cov(\cN)$
many continuous images of $S^{k-1}$.
\end{lem}
\begin{proof}
Let $\kappa = \cov(\cN)=\cof(\cN)$.
Let $\{N^{(k)}_\alpha : \alpha<\kappa\}$ be a
cofinal family of null sets in $\R^k$, $k\in\N$.

For $J\sbst\R^k$, $x\in\R$, and $i<k$, define
$J_{(x,i)} = \{v_1\frown v_2 : \ell(v_1)=i,\ v_1\frown x\frown v_2\in J\}$.
By the Fubini Theorem, for each null set $N\sbst\R^k$ and $i<k$,
$$\tilde N = \{x : (\E i<k)\ N_{(x,i)}\mbox{ is not null in }\R^{k-1}\}\in\cN.$$

We make an inductive construction on $\alpha<\kappa$ of elements $x_\alpha\in\R$
with auxiliary collections $\cF_\alpha$ of null sets, as follows.
For $\alpha<\kappa$ let $\cP_\alpha = \{N^{(k)}_\alpha : k\in\N\}$.
At step $\alpha$ do the following:
\be
\i Choose
$x_\alpha\nin\Union_{\beta<\alpha}\(\Union_{N\in(\cP_\beta\cup\cF_\beta)\sm P(\R)}\tilde N\cup
\Union_{N\in(\cP_\beta\cup\cF_\beta)\cap P(\R)}N\).$
\i Set $\cF_\alpha = \{N_{(x_\alpha,i)} : \beta<\alpha,\ N\in(\cP_\beta\cup\cF_\beta)\sm P(\R),\ i\in\N\}$.
\ee
This is possible because $x_\alpha$ is required to avoid membership in
a union of less than $\cov(\cN)$ many null sets.

Take $S = \{x_\alpha : \alpha<\kappa\}$. Then $S$ is a $\kappa$-Sierpi\'nski set.
Fix $k$. For each null $N\sbst\R^k$, there exists
$\beta<\kappa$ with $N\sbst N^{(k)}_\beta$.
Whenever $\beta<\alpha_0<\dots<\alpha_{k-1}$, and $\pi$ is a permutation on $\{0,\dots,k-1\}$,
$(N^{(k)}_\beta)_{(x_{\alpha_0},\pi\inv(0))} \in\cF_{\alpha_0}$,
thus
$(N^{(k)}_\beta)_{(x_{\alpha_0},\pi\inv(0)),(x_{\alpha_1},\pi\inv(1))} \in\cF_{\alpha_1}$,
\dots,
$(N^{(k)}_\beta)_{(x_{\alpha_0},\pi\inv(0)),\dots,(x_{\alpha_{k-2}},\pi\inv(k-2))} \in\cF_{\alpha_{k-2}}$,
thus $x_{\alpha_{k-1}}\nin (N^{(k)}_\beta)_{(x_{\alpha_0},\pi\inv(0)),\dots,(x_{\alpha_{k-2}},\pi\inv(k-2))}$,
that is, $(x_{\alpha_{\pi(0)}},\dots,x_{\alpha_{\pi(k-1)}})\nin N^{(k)}_\beta$.

Consequently, $S^k\cap N$ is contained in the union of all sets of the form
$S^i\x\{x_\xi : \xi\le\beta\}\x S^{k-i-1}$ ($i< k$)---a union of $|\beta|<\kappa$
copies of $S^{k-1}$---and $\{v\in S^k : v_i=v_j\}$ ($i<j<k$),
which are continuous images of $S^{k-1}$.
\end{proof}
Now assume that $\cov(\cN)=\cof(\cN)=\b$, and let $S$ be a
$\b$-Sierpi\'nski set as in Lemma \ref{sier}.
We will show by induction that for each $k$, $S^k$ satisfies $\sone(\BG,\BG)$.
By \cite{CBC} it is enough to show that for each null $N\sbst\R^k$,
$S^k\cap N$ satisfies $\sone(\BG,\BG)$.
By Lemma \ref{sier} and the induction hypothesis, $S^k\cap N$ is contained in a union of less than
$\b$ many sets satisfying $\sone(\BG,\BG)$.
In \cite{AddQuad} it is shown that
$\sone(\BG,\BG)$ is preserved under taking unions of size
less than $\b$, and in \cite{ideals} it is shown that
$\sone(\BG,\BG)$ is preserved under taking subsets.
This proves the assertion.

So all finite powers of $S$ satisfy $\sone(\BG,\BG)$,
and $G=\<S\>$ works, since as we mentioned
before, $\sone(\BG,\BG)$ is preserved under taking countable
unions.
Finally, by \cite{CBC} $\sone(\BG,\B)$ in all finite powers implies
$\sfin(\BO,\BO)$.
\end{proof}

It seems that the following was not known before.

\begin{cor}
Assume that $\cov(\cN)=\cof(\cN)=\b$. Then
there exists a set of reals satisfying
$\sone(\BG,\BG)$ and $\sfin(\BO,\BO)$, but not $\sone(\O,\O)$.
\end{cor}

\end{document}